\documentclass[graybox]{svmult}

\usepackage{mathptmx}       
\usepackage{helvet}         
\usepackage{courier}        
\usepackage{makeidx}         
\usepackage{graphicx}        
\usepackage{multicol}        
\usepackage[bottom]{footmisc}
\usepackage{amssymb}         

\makeindex             
\usepackage{array}
\usepackage{subfigure} 
\usepackage{enumerate}

\begin{document}

\title*{Mind the Gap: A Study in Global Development through Persistent Homology}
\author{Andrew Banman and Lori Ziegelmeier}
\institute{Andrew Banman \at University of Minnesota, 3 Morrill Hall 100 Church St. S.E., Minneapolis, MN 55455, \email{banma001@umn.edu}
\and Lori Ziegelmeier \at Macalester College, 1600 Grand Avenue, Saint Paul, MN 55105, \email{lziegel1@macalester.edu}}
%
%
\maketitle

\abstract{The Gapminder project set out to use statistics to dispel simplistic notions about global development. In the same spirit, we use persistent homology, a technique from computational algebraic topology, to explore the relationship between country development and geography. For each country, four indicators, gross domestic product per capita; average life expectancy; infant mortality; and gross national income per capita, were used to quantify the development. Two analyses were performed. The first considers clusters of the countries based on these indicators, and the second uncovers cycles in the data when combined with geographic border structure. Our analysis is a multi-scale approach that reveals similarities and connections among countries at a variety of levels. We discover localized development patterns that are invisible in standard statistical methods.
}

\section{Introduction}
\label{sec:intro}

\begin{figure}[t]
\centering
\includegraphics[scale=.35]{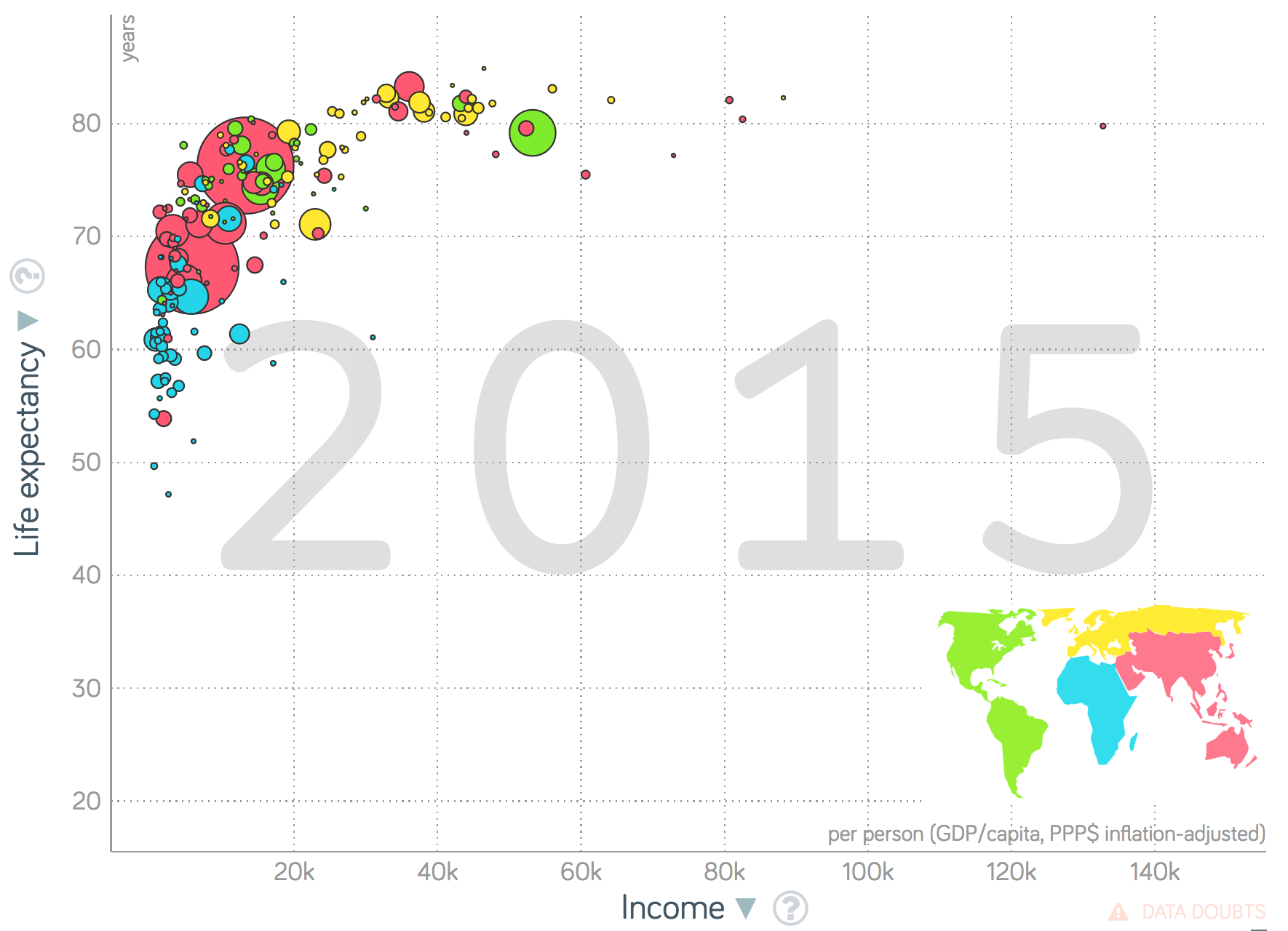}
\caption{``Health \& Wealth of Nations'' chart from the Gapminder World project \cite{GapminderWorld}.}
\label{gapminder-map}
\end{figure}

The Gapminder World \cite{GapminderWorld} project provides a viewpoint of global development through a statistical lens. The first chart that loads in Gapminder plots each country's gross domestic product (GDP) against the life expectancy of its citizens, see Fig.~\ref{gapminder-map}. The project equates GDP per capita with a nation's wealth and life expectancy with its health. Countries are color-coded by their broad geographic region: the Americas, Eurasia, etc. A time lapse animation shows countries transitioning along a common trajectory towards more health and wealth, telling a common story about global development. However, it is not clear what role geography plays in this trend. While one may say that most African nations lag behind most Eurasian nations, it is difficult to draw any finer conclusions solely from these two statistics, as each region spans a large range of the development statistics. Furthermore, Gapminder's pre-determined regions have been chosen according to a convention rather than from the data. For instance, it splits the African continent into Northern and Sub-Saharan regions, isolates India and a few of its neighbors, and joins Australia with Southern Asian countries. These regions do not necessarily align with regions of differing development.

We seek a quantifiable, fine-grained, and unbiased method to analyze development and geographic trends in this data. Persistent homology \cite{barcodes,carlsson2009topology, edelsbrunner2008persistent} gives us tools to uncover the structure of high-dimensional, complicated data, revealing groups (connected components) and cycles (loops) in the data at multiple scales. Persistent homology has been used to understand the topological structure of data arising from applications including computer vision, biological aggregations, brain structure, among many others \cite{imagewebs, windowsandpersistence, hippocampalPH, corticalsurfacePD, visionTDA, Swarms}. In particular, the paper \cite{DBLP:journals/corr/StolzHP16} analyzes data related to the recent, so-called, ``Brexit" referendum using persistent homology.

We use persistent homology to expand on Gapminder's study of health and wealth statistics. We explore two methods (1) computing the connected components of the indicators of GDP per capita and life expectancy as well as infant mortality and gross national income per capita and (2) adding the underlying geography to the indicators by constructing a weighted graph based on country borders to observe cycles in the data. The structure of the data is uncovered at multiple scales. Our analyses reveal that there are connections among countries at a variety of levels and show subtleties with country similarities and differences, as well as loops formed by countries geographically linked. This provides a more nuanced view than simply the ``first'' versus ``third'' world paradigm, a construction that divides the world into discrete sets of developed and undeveloped countries \cite{nationsonline}.

The remainder of this paper proceeds as follows. Background on the computational approach of persistent homology is discussed in Section \ref{sec:background}. Section \ref{methods} outlines the indicators we use to quantify health and wealth of nations, and our implementation of persistent homology on these indicators. We analyze the results of these computations in Section \ref{sec:barcodes}. Conclusions and future work are discussed in Section \ref{sec:conclusion}.

\section{Background on Persistent Homology} \label{sec:background}

\emph{Persistent homology} is a computational approach to topology that encodes a parameterized family of homological features such as connected components, loops, trapped volumes, etc of a topological space. It allows one to answer basic questions about the structure of point clouds at multiple scales. As such, it can uncover the ``shape" of data. Broadly, this procedure involves (1) interpreting a point cloud as a noisy sampling of a topological space, (2) creating a global object by forming connections between proximate points based on a scale parameter, (3) determining the topological structure made by these connections, and (4) looking for structures that persist across different scales. For foundational material and overviews of computational homology in the setting of persistence, see \cite{edelsbrunner2008persistent,Edelsbrunner10,barcodes,carlsson2009topology,computingPH}.

Beginning with a finite set of data points, a nested sequence of simplicial complexes indexed by a parameter $\epsilon$ may be created by taking the vertices as the data points and forming a $k$-simplex whenever $k+1$ points are pairwise within distance $\epsilon$. This procedure is known as the \emph{Vietoris-Rips} (VR) complex which is often used for its computational tractability \cite{barcodes}. Fixing a field $\mathbb F$, one builds a chain complex of vector spaces over $\mathbb F$ for each simplicial complex. For each pair $\epsilon_1<\epsilon_2$, there is a pair of simplicial complexes, $S_{\epsilon_1}$ and $S_{\epsilon_2}$, and an inclusion map $j:S_{\epsilon_1} \hookrightarrow S_{\epsilon_2}$. This inclusion map induces a chain map between the associated chain complexes which further induces a linear map between the corresponding $k^{th}$ homology vector spaces.  The dimension of the $k^{th}$ homology vector space is known as the $k^{th}$ \emph{Betti number} $\beta_k$ and corresponds to the number of connected components, loops, trapped volumes, etc. of a simplicial complex for $k=0,1,2,\ldots$, respectively.

The $k^{th}$ barcode is a way of presenting Betti numbers across multiple scales $\epsilon$ \cite{barcodes}. From the barcode, one can visualize the number of independent homology classes that persist across a given filtration interval $[\epsilon_b,\epsilon_d]$ as a function of the scale $\epsilon$. See the top row of Fig.~\ref{fig:CountryBorderBarcodes2d} for an example $\beta_0$ barcode and the bottom row of Fig.~\ref{fig:CountryBorderBarcodes2d} for an example $\beta_1$ barcode. Each horizontal bar begins at the scale where a topological feature first appears (``is born") and ends at the scale where the feature no longer remains (``dies"). The $k^{th}$ Betti number at any given parameter value $\epsilon$ is the number of bars that intersect the vertical line through $\epsilon$.  For $\beta_0$ in our setting, there will be a distinct bar for each data point at small values of $\epsilon$, as the simplicial complex $S_\epsilon$ consists only of isolated points. At large values of $\epsilon$, only one bar remains as all data will eventually connect into a single component.

The idea of persistence is to not only consider the homology for a single specified choice of parameter $\epsilon$ but rather, track topological features through a range of parameters.  Those which persist over a large range of values are considered signals of underlying topology, while the short lived features are taken to be noise inherent in approximating a topological space with a finite sample \cite{carlsson2009topology}.

\section{Methods}
\label{methods}

There are many ways to quantify the health and wealth of nations. We study four development indicators: gross domestic product (GDP) per capita\footnote{Gross Domestic Product per capita by Purchasing Power Parities (in international dollars, fixed 2011 prices). The inflation and differences in the cost of living between countries has been taken into account \cite{worldbankGNI}.}, life expectancy\footnote{The average number of years a newborn child would live if current mortality patterns were to stay the same \cite{worldbankGDP}.}, rate of infant mortality\footnote{The probability that a child born in a specific year will die before reaching the age of one, if subject to current age-specific mortality rates. Expressed as a rate per 1,000 live births \cite{gbd2013}.}, and gross national income (GNI) per capita\footnote{Gross national income converted to international dollars using purchasing power parity rates \cite{UNICEF}.}. These indicators were chosen because (1) we believe them to be broad indicators of health and wealth, and (2) recent data is available for a large set of countries in each indicator.

We consider this data in two sets: what we will call the four-dimensional ($\mathbb{R}^4$) data comprising all four indicators and the two-dimensional ($\mathbb{R}^2$) data comprising only GDP/capita and life expectancy. The raw $\mathbb{R}^2$ data---before scaling as discussed below---generates the Gapminder chart, see Fig.~\ref{gapminder-map}, allowing a comparison of our results to the chart.

The frequency of reporting and currency of statistics can vary dramatically by country so any result necessarily carries the ``according to available data'' qualifier. We construct our data sets by taking the most recent value for each indicator corresponding to a country\footnote{Most data comes from years 2015, 2016, with others as early as 2005. See Table~\ref{data_years} in Appendix A.}. Countries with no available data for one or more indicators in this time frame are excluded from the data set. This yields data comprising 194 countries in the $\mathbb{R}^2$ set and 179 countries in $\mathbb{R}^4$. See Table~\ref{indicatorStats} for statistics such as the maximum, minimum, median, mean, and standard deviation for the raw data of the indicators.

\begin{table}[h!]
 \caption{Statistics of each indicator: GDP per capita (GDP), Life Expectancy (LE), Infant Mortality rate (IM), and GNI per capita (GNI). The first five statistics correspond to the raw data; the last corresponds to the attenuated and scaled data. Naturally, high GDP, GNI, and life expectancy are favorable, whereas high infant mortality rate is unfavorable.}
 \centering
 \begin{tabular}{lcccccc}
 \hline\noalign{\smallskip}
 Indicator & Max & Min & Median& Mean & Stand Dev & Scaled Mean \\
 \noalign{\smallskip}\hline\noalign{\smallskip}
GDP & 148374 & 599 & 11903 & 18972 & 21523 & -0.476 \\
LE & 84.8 & 48.86 & 74.5& 72.56 & 7.74 & 0.296 \\
IM & 96 & 1.5 & 23.89 & 15 & 21.9 & 0.528 \\
GNI & 87030 & 350 & 8360 & 13596 & 15399 & -0.431\\
 \noalign{\smallskip}\hline\noalign{\smallskip}
 \end{tabular}
 \label{indicatorStats}
\end{table}

We consider the relative health and wealth of countries, and the presence of extreme outliers in GDP obscures this relationship. Rather than exclude these countries outright, we modulate their values to two standard deviations from the mean. Alternatively, we could have taken the logarithm of GDP to bring the outliers closer to the bulk. However, this option has the undesirable consequence of exaggerating the distance between countries with very low GDP and understating the distance between higher GDP countries. For our purposes, it made more sense to collect the richest countries into one group at the extreme of the spectrum and likewise for the poorest. The same attenuation was done for the GNI per capita indicator.

Each indicator is then re-scaled to $[-1,1]$. The range $[-1,1]$ was chosen to give a normative representation of each indicator, in which -1 is least favorable and 1 is most favorable, e.g. the country with lowest life expectancy has -1 for that dimension and the country with lowest infant mortality has 1 in that dimension. Note that this does not imply zero is the average value for any indicator. There are many more relatively low GDP countries, even after attenuating outliers, see Table~\ref{indicatorStats}. This scaling is required to ensure each indicator carries equal weight in the persistent homology calculations. Otherwise GDP/capita and GNI/capita would completely obscure any features in life expectancy and infant mortality rate because they are orders of magnitude larger in conventional units.

For our calculations, we use the TDA library in R \cite{TDApackage}. This library provides an API to create a filtered simplicial complex upon which to calculate the persistent homology. The final result of the computation is a list of persistence intervals $[\epsilon_b,\epsilon_d]$, neatly displayed in a homology barcode, where each interval indicates a homological feature that is born at $\epsilon_b$ and dies at $\epsilon_d$. In this section, we outline our procedure for computing persistent homology of our data. In the next, we analyze the results.

\begin{figure}
\centering
\subfigure[][]{
\label{fig:CountryBarcodes2d}
\includegraphics[width=0.49\linewidth]{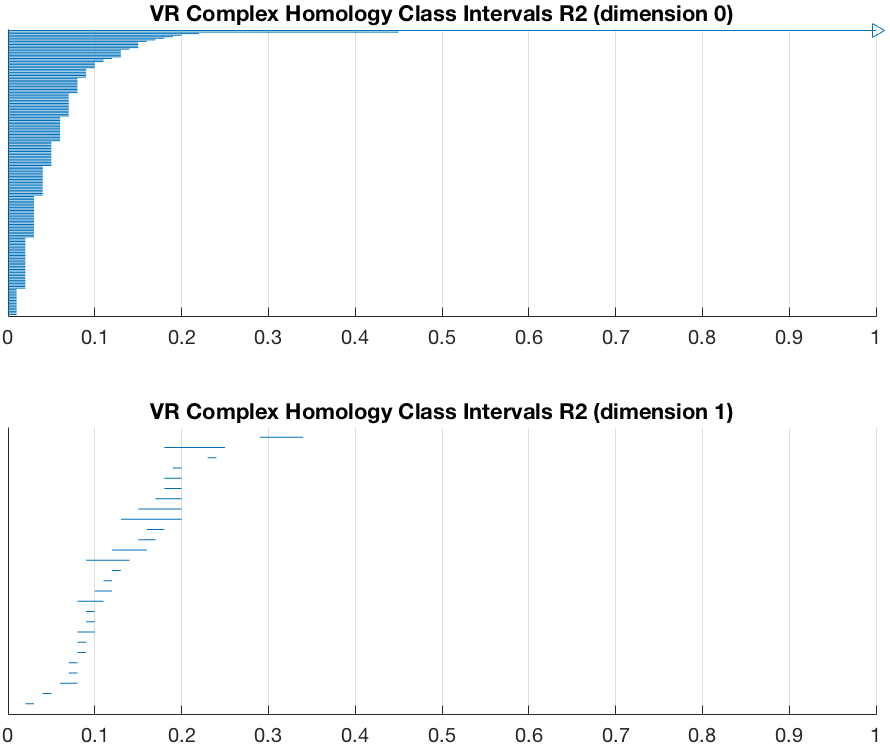}}
\subfigure[][]{
\label{fig:CountryBarcodes4d}
\includegraphics[width=0.49\linewidth]{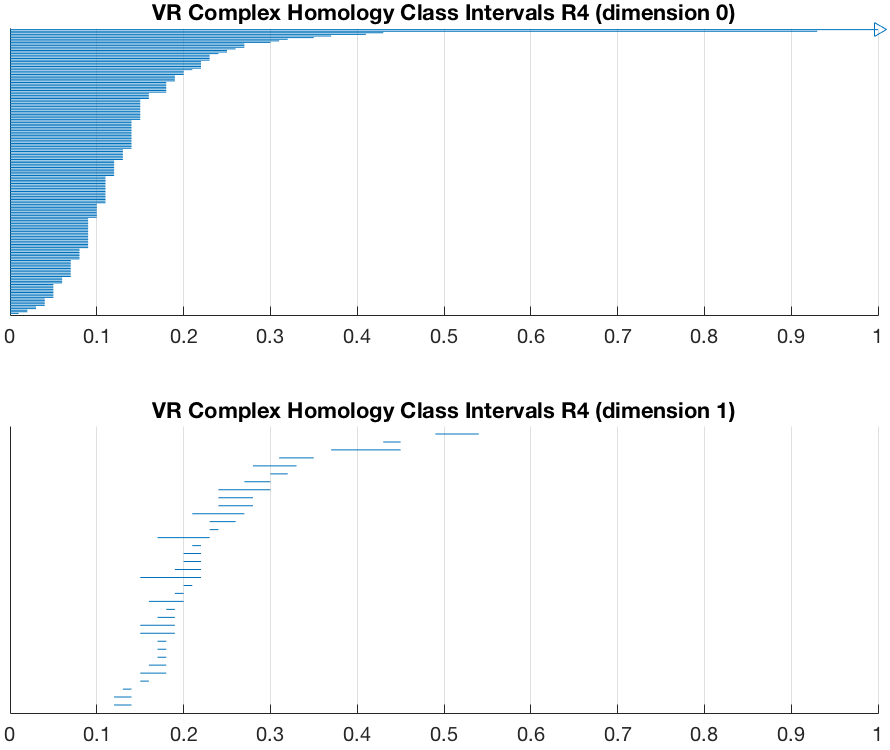}} \\
\caption[]{Zero-order (top) and first-order (bottom) persistent homology barcodes of the VR complex stream over point cloud in $\mathbb{R}^{|I|}$:
\subref{fig:CountryBarcodes2d} $I=\{\textrm{GDP}, \textrm{LifeExp}\}$ and
\subref{fig:CountryBarcodes4d} $I=\{\textrm{GDP}, \textrm{LifeExp}, \textrm{InfMort}, \textrm{GNI}\}$.}
\label{fig:CountryBarcodes}
\end{figure}

For our first experiment, we interpret each set of countries as a point cloud with each indicator value as a dimension. We then apply the Euclidean metric to define the distance between two countries $x$ and $y$ over a set of indicators $I$:
$$ d_I:\mathbb{R}^{|I|} \rightarrow \mathbb{R} $$
$$ d_I(x,y) = \sqrt{\sum\limits_{i\in I} (x_i - y_i)^{2}} $$
We use TDA to construct a stream of VR complexes from these point clouds over a range of filtration values $\epsilon\in[0,1.0]$. Fig.~\ref{fig:CountryBarcodes} shows the zero-order and first-order barcodes of the VR streams for the two sets of indicators ($\mathbb{R}^2$ on the left and $\mathbb{R}^4$ on the right).

\begin{figure}
\centering
\subfigure[][]{
\label{fig:CountryBorderBarcodes2d}\includegraphics[width=0.49\linewidth]{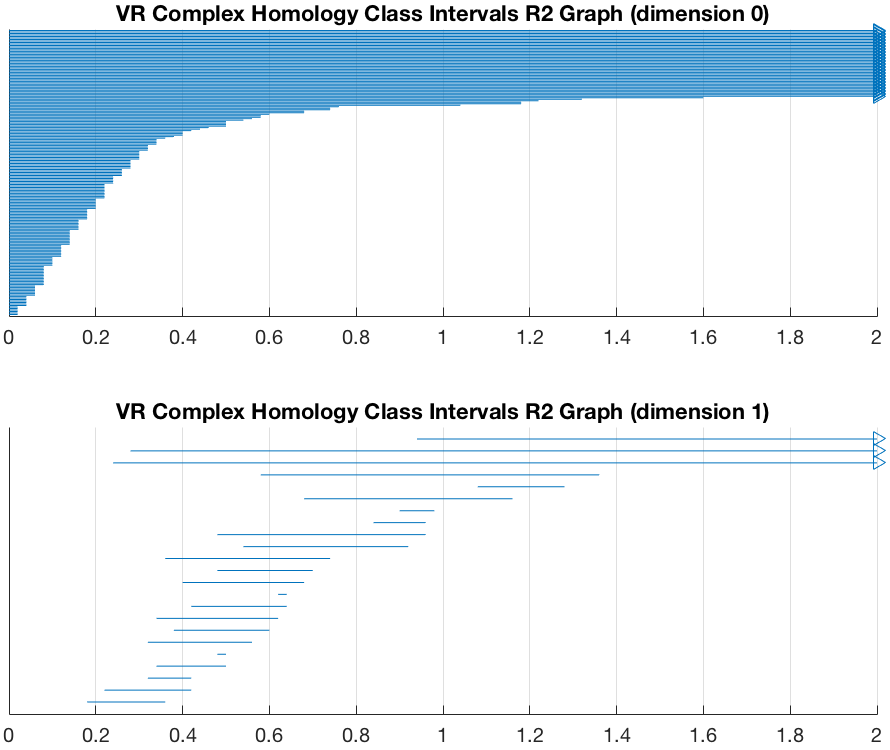}}
\subfigure[][]{
\label{fig:CountryBorderBarcodes4d}
\includegraphics[width=0.49\linewidth]{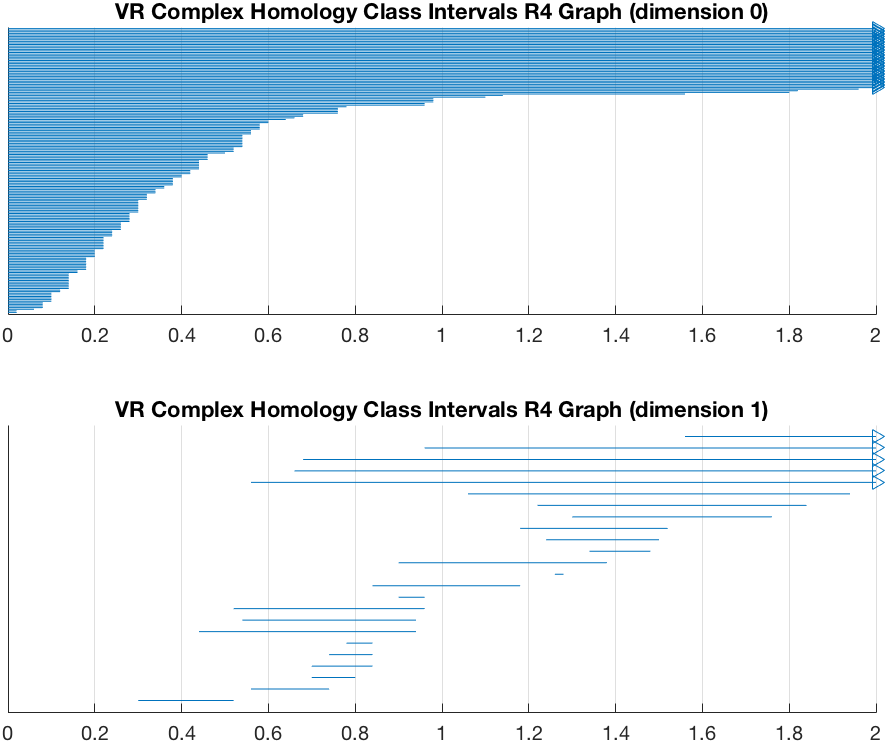}} \\
\caption[]{Zero-order (top) and first-order (bottom) persistent homology barcodes of the VR complex stream over country border graph with distance $d_I$ as the edge weight.
\subref{fig:CountryBorderBarcodes2d} $I=\{\textrm{GDP}, \textrm{LifeExp}\}$ and
\subref{fig:CountryBorderBarcodes4d} $I=\{\textrm{GDP}, \textrm{LifeExp}, \textrm{InfMort}, \textrm{GNI}\}$.
}
\label{fig:CountryBorderBarcodes}
\end{figure}

For the second experiment, we add the geographic structure to the data by constructing a weighted graph over the countries and their borders. From country border data \cite{geonames}, we define an adjacency matrix $A$
\[ A_{(i,j)} = \left\{
\begin{array}{ll}
1 & {\rm if\ countries\ }{i,j\rm\ share\ a\ border}, \\
0 & {\rm if\ countries\ }{i,j\rm\ do\ not\ share\ a\ border} \end{array} \right.\]
from which we arrive at the distance matrix $D$ for a set of indicators $I$,
\[ D_{I(i,j)} = \left\{
\begin{array}{ll}
d_I(i,j) & {\rm if\ }{A_{i,j}=1}, \\
\infty & {\rm if\ }{A_{i,j}=0} \end{array} \right.\]
where, for practicality, infinity is set to be a number larger than the maximum filtration value. This maximum is chosen to be large enough to display the entire set of intervals. We then compute the persistent homology of the explicit metric space defined by $D_I$.\footnote{It has been observed that, for the VR complex, the metric in question need not actually be a metric as it is not a requirement to satisfy the triangle inequality \cite{JavaPlexTutorial}. The construction described here is also known as a weighted rank clique complex. For example, see \cite{DBLP:journals/corr/StolzHP16}.} The zero-order and first-order persistent homology barcodes for the weighted graphs over $\mathbb{R}^2$ data and $\mathbb{R}^4$ data are shown in Fig.~\ref{fig:CountryBorderBarcodes}. In this framework incorporating geographic structure, our focus is on the first-order features.

Generally, longer intervals are construed to represent more significant homology classes while short intervals are noise in the data. Statistically significant intervals can be quantitatively determined by the methods presented in \cite{fasy2014}. However, we shall see that even relatively short intervals in the first-order barcode reveal interesting patterns in the development indicators. On the other hand, intervals in Fig.~\ref{fig:CountryBorderBarcodes} that persist through the full range of the filtration are less interesting to us as they relate to the inherent border graph structure. These ``infinite'' intervals in the dimension-0 barcode indicate island nations that share no borders with other countries. Since their distance to all other countries is infinite, they remain distinct components in the VR complex. The infinite intervals in the dimension-1 barcodes indicate homology classes inherent to the country border graph. The three infinite intervals in Fig.~\ref{fig:CountryBorderBarcodes2d} identify the Black, Caspian, and Mediterranean seas. Fig.~\ref{fig:CountryBorderBarcodes4d} has two additional intervals that exist because two countries (South Sudan and Zimbabwe) were dropped from the data set as not all four indicators were present, creating holes in the graph not unlike an inland sea. That these features are identified is a good sanity check for the method.

\section{Parsing the Barcodes}
\label{sec:barcodes}

\subsection{Clustering of Development Groups} \label{sec:H0clustering}

Zero-order persistent homology can be viewed as a clustering algorithm, where the connected components of a simplicial complex represent clusters in the data. In fact, these components are equivalent to the clusters of the hierarchical method of single-linkage clustering. In Fig.~\ref{gapminder-map}, we see a clustering chosen by Gapminder. In this section, we describe the clusters found using zero-order persistent homology present in the barcode of Fig.~\ref{fig:CountryBarcodes}, focusing on the first experiment which only relies on distance between indicators and does not incorporate the country border information. In Appendix B, we present clusters selected by the classic $K$-means algorithm. Each of these methods results in different clusters. However, we observe that viewing clusters at multiple scales and adding more indicators provides additional insight into relations among countries in terms of health and wealth.

We examine the clusters found using dimension-0 persistent homology by extracting the elements in each component of the simplicial complex for a particular filtration value, see the top row of Fig.~\ref{fig:CountryBarcodes}. One may imagine drawing a vertical slice through the dimension-0 barcode at a given $\epsilon$ to select the components. We then extract the list of countries comprising each component using a union-find algorithm. The Betti number can be viewed as a function of the filtration value, $\beta_k(\epsilon)$. When $\epsilon = 0$, each country is an isolated point, and hence, $\beta_0(0) = 194$ for the $\mathbb{R}^2$ data and $\beta_0(0) = 179$ for the $\mathbb{R}^4$ data. All countries in the point cloud eventually merge into a single connected component. This occurs at approximately $\epsilon=0.45$ for the $\mathbb{R}^2$ data and $\epsilon=0.92$ for the $\mathbb{R}^4$ data, as seen in the barcodes where only one bar remains.

\begin{figure}%
\centering
\subfigure[][]{%
\label{2dMaps-a}%
\includegraphics[width=0.48\linewidth]{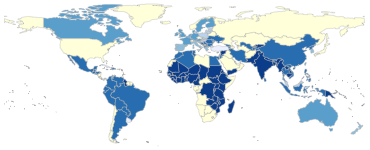}}%
\hspace{8pt}%
\subfigure[][]{%
\label{2dMaps-b}%
\includegraphics[width=0.48\linewidth]{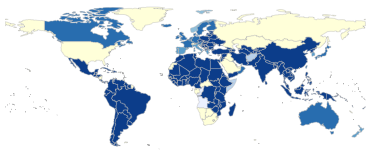}} \\
\subfigure[][]{%
\label{2dMaps-c}%
\includegraphics[width=0.48\linewidth]{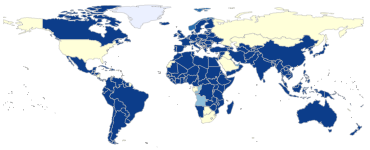}}%
\hspace{8pt}%
\subfigure[][]{%
\label{2dMaps-d}%
\includegraphics[width=0.48\linewidth]{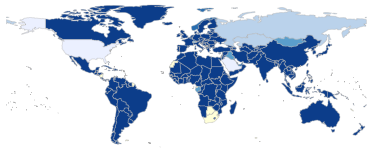}}%
\caption[]{World map depicting clusters found using dimension-0 persistent homology of the VR complex of the $\mathbb{R}^2$ data at various filtrations. The six largest connected components are displayed in shades of blue (darker indicates larger) while other countries not in these clusters are displayed in yellow.
\subref{2dMaps-a} $\epsilon=0.08$ with six largest components consisting of 54, 52, 14, 10, 8 countries among 41 total distinct clusters;
\subref{2dMaps-b} $\epsilon=0.10$ with six largest components consisting of 132, 18, 10, 6, 2 countries among 25 total distinct clusters;
\subref{2dMaps-c} $\epsilon=0.12$ with six largest components consisting of 164, 6, 2, 2, 2 countries among 19 total distinct clusters;
\subref{2dMaps-d} $\epsilon=0.14$ with six largest components consisting of 170, 8, 3, 2, 2 countries among 13 total distinct clusters.}
\label{2dMaps}
\end{figure}

\begin{figure}
\centering
\subfigure[][]{
\label{4dMaps-a}
\includegraphics[width=0.4\linewidth]{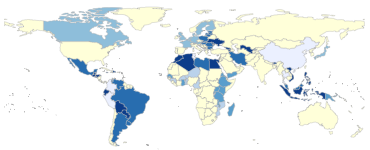}}
\hspace{8pt}
\subfigure[][]{
\label{4dMaps-b}
\includegraphics[width=0.47\linewidth]{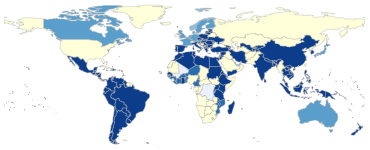}} \\
\subfigure[][]{
\label{4dMaps-c}
\includegraphics[width=0.47\linewidth]{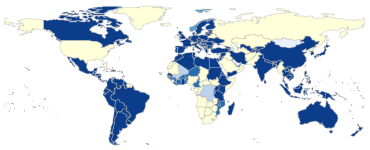}}
\hspace{8pt}
\subfigure[][]{
\label{4dMaps-d}
\includegraphics[width=0.47\linewidth]{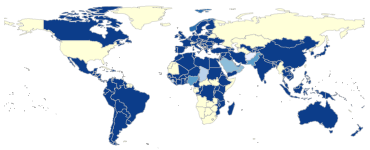}}
\caption[]{World map depicting clusters found using dimension-0 persistent homology of the VR complex of the $\mathbb{R}^4$ data at various filtrations. The six largest connected components are displayed in shades of blue (darker indicates larger) while other countries not in these clusters are displayed in yellow.
\subref{4dMaps-a} $\epsilon=0.14$ with six largest components consisting of 26, 18, 13, 11, 8 countries among 76 total distinct clusters;
\subref{4dMaps-b} $\epsilon=0.16$ with six largest components consisting of 101, 13, 12, 7, 3 countries among 44 total distinct clusters;
\subref{4dMaps-c} $\epsilon=0.18$ with six largest components consisting of 115, 13, 7, 3, 2 countries among 40 total distinct clusters;
\subref{4dMaps-d} $\epsilon=0.20$ with six largest components consisting of 133, 7, 5, 3, 2 countries among 33 total distinct clusters.}
\label{4dMaps}
\end{figure}

Fig.~\ref{2dMaps} and Fig.~\ref{4dMaps} display the six\footnote{The choice of six is to coincide with the six clusters in the Gapminder project, see Fig. \ref{gapminder-map}.} components that contain the largest number of countries in the cluster at a variety of filtration scales for the $\mathbb{R}^2$ and $\mathbb{R}^4$ data, respectively. We further inspect these components in detail below.

First, we consider the large-scale structure of the data. For the $\mathbb{R}^2$ point cloud there are 170 countries in a single connected component at $\epsilon=0.14$, eight countries in the next largest, and the remaining countries isolated in small components. We may say this large cluster is the dominant feature of the data. The $\mathbb{R}^4$ point cloud shows the same behavior. Fig.~\ref{2dMaps} and Fig.~\ref{4dMaps} show how quickly this dominant component grows at early filtration values. At no point do we observe two dominant clusters capturing a combined majority of countries.

Thus, the dimension-0 clustering shows that countries of the world may not be neatly divided into ``first world'' and ``third world'' categories with this method.\footnote{The clustering presented in Appendix B results in different clusters, which more closely align with this simplistic notion.} The vast majority of countries are statistically quite similar to another country, which itself is similar to some other country, and so on. The result is a gradient in health and wealth statistics, rather than a discrete grouping. This is easily visualized in the Gapminder chart Fig. \ref{gapminder-map}. One sees the countries of the world arrayed along a gradient from poorer countries with less longevity to richer, longer living countries. Persistent homology clustering captures this gradient as the resulting clusters from this method connect points to their nearest neighbors which each connect to their nearest neighbors and so on. This may result in long clusters whose elements at the ends of a cluster may be quite different from one another but are connected through their neighbors.

We also examine the small-scale structure by looking at smaller $\epsilon$ cross-sections. Fig.~\ref{2dMaps} and Fig.~\ref{4dMaps} show a sampling of clusters for early filtration values, before most countries are joined up into one dominant cluster. Consider the clusters in $\mathbb{R}^2$ at $\epsilon=0.08$, shown in Fig.~\ref{2dMaps-a} and detailed in Table~\ref{clusterTable}. While most countries fall into connected components of one to four countries, there are six larger components that capture 138 countries. Because these clusters only exist at a small scale, the countries in each cluster must be quite close in the data. Hence, we may conceive of these groups as sets of very similar countries according to the indicators. This clustering makes a distinction between groups of countries with varying GDP/capita and similar life expectancy. Observe clusters 2-4 have similar life expectancy but a wide range of increasing GDP. Likewise, clusters 5, 6 have almost the same LE but a 0.4 gap in GDP. From this result we may conclude there is nuance in development among poor countries that may be obfuscated by the "third-world" identifier.

\begin{table}
\scriptsize
\caption{Countries comprising the largest connected components in the VR complex at filtration $\epsilon=0.08$ over $\mathbb{R}^2$ and the corresponding means of scaled indicators, GDP/capita (GDP) and life expectancy (LE), for each cluster. Clusters are listed in ascending GDP order, for clarity in comparison.}
\label{clusterTable}
\begin{tabular}{p{8cm}rr}
\hline\noalign{\smallskip}
Countries (ISO2) & GDP & LE \\
\noalign{\smallskip}\svhline\noalign{\smallskip}
Bangladesh, Kyrgyzstan, Cambodia, Mauritania, Micronesia Fed. Sts., Nepal, Syria, Gambia, Comoros, Myanmar, Sudan, Sao Tome and Principe, India, Laos, Marshall Islands, Guyana, Pakistan, Ghana, Nigeria, Yemen Rep., Djibouti, Kenya, Senegal, Tanzania, Vanuatu, Haiti, Liberia, Madagascar, Solomon Islands, Ethiopia, Rwanda, Benin, Kiribati, Burkina Faso, Burundi, Congo Dem. Rep., Niger, Papua New Guinea, Togo, Uganda, Zimbabwe, Eritrea, Mali, Malawi, Guinea, Cote d'Ivoire, Cameroon, Sierra Leone, Mozambique, Chad, Zambia, South Sudan, Guinea-Bissau, Fiji & -0.93 & -0.15 \\
Albania, Bosnia and Herzegovina, Colombia, Jordan, Sri Lanka, Tunisia, Peru, Macedonia FYR, Barbados, China, Dominican Rep., Algeria, Ecuador, Montenegro, Serbia, Thailand, Bulgaria, Brazil, Iran, Venezuela, Mauritius, Mexico, Romania, Argentina, Saint Lucia, Armenia, Jamaica, Paraguay, El Salvador, Morocco, Vietnam, Bolivia, Bhutan, Cape Verde, Georgia, Guatemala, Honduras, Moldova, Samoa, Belize, Ukraine, Indonesia, Philippines, Saint Vincent and the Grenadines, Egypt, Grenada, Tonga, Uzbekistan, Tajikistan, Korea Dem. Rep., Timor-Leste, Palestine & -0.69 & 0.44 \\
Antigua and Barbuda, Croatia, Uruguay, Cuba, Panama, Turkey, Lebanon & -0.37 & 0.63 \\
Estonia, Poland, Slovak Republic, Hungary, Latvia, Malaysia, Lithuania, Seychelles & -0.19 & 0.53 \\
Cyprus, Malta, Slovenia, Israel, Spain, Italy, Korea Rep., New Zealand, Portugal, Greece & -0.02 & 0.83 \\
Austria, Australia, Canada, Germany, Denmark, Netherlands, Sweden, Belgium, Taiwan, Finland, France, United Kingdom, Bahrain, Ireland & 0.38 & 0.80 \\
\noalign{\smallskip}\hline\noalign{\smallskip}
\end{tabular}
\end{table}

One advantage of persistent homology as a clustering algorithm is the total lack of bias in the origination of each cluster. Further, a smaller filtration $\epsilon$ yields a finer clustering, whereas a relatively large $\epsilon$ reveals a coarser structure of the data allowing for a multi-scale analysis. However, the algorithm is highly sensitive to ``bridge" structures that connect one cluster to another, destroying distinct components. A bridge in our data might be a relatively poor country with high longevity connecting to a relatively wealthy country with similar longevity, thus joining a cluster of poor countries with a cluster of wealthier countries.

\subsection{Local Development Patterns}

First-order homology classes represent cycles in the data, often visualized as loops around a hole in the point cloud or graph. The dimension-1 barcode intervals tell us over what range of filtration values these cycles exist. Software provides a list of generating simplicies for each homology class, which we parse as a list of generating countries. These generators tell us where in the world the cycle exists. \footnote{The generating countries are not guaranteed to be minimal in a geometric sense; they can make up any loop through the connected component that contains the homology class. One can find the minimal loop by examining the weight of its internal edges.}

\begin{figure}[t]
\centering
\includegraphics[scale=.70]{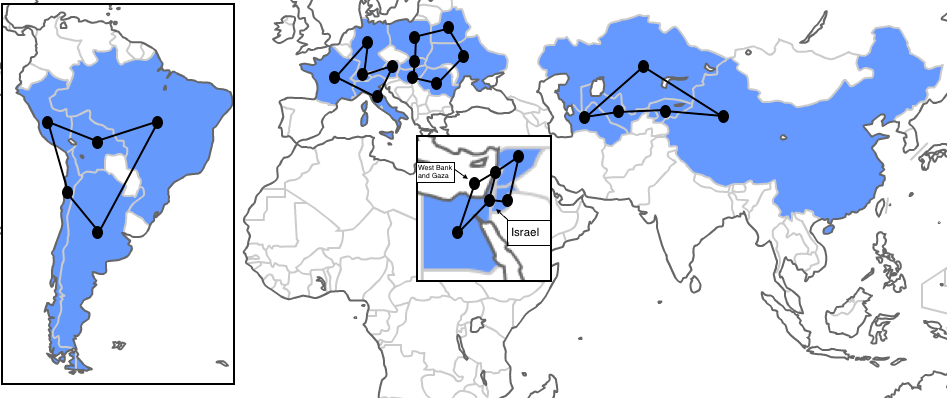}
\caption{Map of six cycles in the country border graph with distance $d_I$, where $I=\{\textrm{GDP}, \textrm{LifeExp}\}$, as the edge weight. Software-provided generators for each cycle are shown, and the involved countries are shaded.}
\label{cyclemap}
\end{figure}

Our main focus for dimension-1 homology is on the weighted border graphs (Fig.~\ref{fig:CountryBorderBarcodes}). The most discernible intervals are those persisting through the full range of the filtration. As discussed in Section~\ref{methods}, these infinite loops describe the topology of the border graph itself. More interesting are the cycles not inherent to the graph structure. These cycles exist because of a pattern of similarity between country neighbors in the indicator data. We map out the generating countries of the six longest-persisting (non-infinite) cycles from the $\mathbb{R}^2$ data in Fig.~\ref{cyclemap} and further highlight the countries generating two of these cycles in Tables~\ref{cycle1} and \ref{cycle2}. The cycles are distinguished by a periodic pattern in the data, in which a ``maximal'' country has the greatest value in one or more indicators; a ``minimal'' country has the least in these indicators; and the connected countries have intermediary values.

\begin{table}
\parbox{.45\linewidth}{
\caption{Generating countries of the South America cycle in the $\mathbb{R}^2$-weighted graph from the dimension-1 barcode interval $[0.34, 0.62)$ in Fig.\ref{fig:CountryBorderBarcodes2d}.}
\begin{tabular}{lrr}
\hline\noalign{\smallskip}
Country & GDP & LE \\
\noalign{\smallskip}\hline\noalign{\smallskip}
Chile       & -0.29 & 0.71 \\
Peru        & -0.63 & 0.72 \\
Bolivia     & -0.81 & 0.37 \\
Brazil      & -0.52 & 0.43 \\
Argentina   & -0.45 & 0.55 \\
\noalign{\smallskip}\hline\noalign{\smallskip}
\end{tabular}
\label{cycle1}
}
\hfill
\parbox{.45\linewidth}{
\caption{Generating countries of the North Africa cycle in the $\mathbb{R}^2$-weighted graph from the dimension-1 barcode interval $[0.85, 0.97)$ in Fig.\ref{fig:CountryBorderBarcodes2d}.}
\begin{tabular}{lrr}
\hline\noalign{\smallskip}
Country & GDP & LE \\
\noalign{\smallskip}\hline\noalign{\smallskip}
Libya       & -0.46 &  0.36 \\
Niger       & -0.99 & -0.31 \\
Mali        & -0.96 & -0.36 \\
Mauritania  & -0.89 &  0.17 \\
Algeria     & -0.58 &  0.54 \\
\noalign{\smallskip}\hline\noalign{\smallskip}
\end{tabular}
\label{cycle2}
}
\end{table}

For example, consider the South American cycle in Table~\ref{cycle1}. A generator of this cycle has Chile with the highest GDP per capita and life expectancy at $(-0.29, 0.71)$ and Bolivia with the lowest at $(-0.81, 0.37)$. Each indicator decreases as you follow the cycle from Chile to Bolivia\footnote{There is a slight deviation from monotonic decrease in life expectancy at Peru. These deviations are not uncommon, but do not detract from the maximal-minimal pattern we observe.}, and increases on the way back around. The same result holds for the North Africa cycle displayed in Table~\ref{cycle2}, where Libya is maximal and Mali is minimal. Thus, persistent homology has identified a set of nearby countries that conform to a cycle in both health and wealth statistics. In other words, we've identified a maximal and a minimal country, local to a connected region, where neighboring countries exist on a gradient between the two poles.

The same patterns are found among the $\mathbb{R}^4$ cycles. Almost all of the cycles found in the $\mathbb{R}^4$ case are present in the $\mathbb{R}^2$ case, albeit with different filtration values for $\epsilon_b$ and $\epsilon_d$. Some cycles not repeated in $\mathbb{R}^4$ are those involving countries that don't exist in the smaller $\mathbb{R}^4$ data set, due to missing data in the added indicators. In some cases the set of generating countries changed, but the core members - i.e. the maximal and minimal countries - are the same.

Additional cycles in the $\mathbb{R}^4$-weighted graph were created or extended by adding the new indicators, infant mortality and GNI/capita (PPP). Table~\ref{cycle3} and Table~\ref{cycle4} show the adjusted indicator values for two such cycles. The first shows periodic patterns in each of the four indicators from maximal Libya to minimal Chad. This cycle length grew from 0.50 in the $\mathbb{R}^2$ case to 0.90 in $\mathbb{R}^4$. The second cycle is an example of a cycle that wasn't already present in the $\mathbb{R}^2$-weighted graph. The four countries are very close in GDP/capita and life expectancy, but there is strong periodic behavior in infant mortality from maximal Senegal to minimal Mali.

\begin{table}
\parbox{.45\linewidth}{
 \caption{Cycle from Libya to Chad found in the country border graph with weight $d_I$, where $I$=\{GDP/capita (GDP), life expectancy (LE), infant mortality (IM), GNI/capita (GNI)\}. Parsed from the interval persisting over $[1.10,1.96)$ in Fig.~\ref{fig:CountryBarcodes4d}.}
 \begin{tabular}{lrrrr}
 \hline\noalign{\smallskip}
 Country & GDP & LE & IM & GNI \\
 \noalign{\smallskip}\hline\noalign{\smallskip}
 Libya  & -0.46 &  0.36 &  0.79 & -0.28 \\
 Sudan  & -0.89 &  0.05 &  0.02 & -0.93 \\
 Chad   & -0.95 & -0.49 & -0.77 & -0.96 \\
 Niger  & -0.99 & -0.31 & -0.18 & -0.98 \\
 \noalign{\smallskip}\hline\noalign{\smallskip}
 \end{tabular}
 \label{cycle3}
}
\hfill
\parbox{.45\linewidth}{
 \caption{Cycle from Senegal to Mali found in the country border graph with weight $d_I$, where $I$=\{GDP/capita (GDP), life expectancy (LE), infant mortality (IM), GNI/capita (GNI)\}. Parsed from the interval persisting over $[0.57,0.75)$ in Fig.~\ref{fig:CountryBarcodes4d}.}
 \begin{tabular}{lrrrr}
 \hline\noalign{\smallskip}
 Country & GDP & LE & IM & GNI \\
 \noalign{\smallskip}\hline\noalign{\smallskip}
 Mauritania  & -0.89 &  0.17 & -0.35 & -0.91 \\
 Senegal    & -0.95 & -0.07 &  0.15 & -0.93 \\
 Guinea     & -0.98 & -0.40 & -0.26 & -0.97 \\
 Mali       & -0.96 & -0.36 & -0.54 & -0.97 \\
 \noalign{\smallskip}\hline\noalign{\smallskip}
 \end{tabular}
 \label{cycle4}
}
\end{table}

The relative scales when cycles are born and die reveal how similar the member countries are to one another; the sooner in the filtration they appear, the more similar we can say they are. This follows from the observation that the beginning of an interval is equivalent to the maximum weight of the cycle's edges. Then, we expect countries in later cycles to be further apart, i.e. less similar in the data, than countries in earlier cycles. The South American cycle in Table~\ref{cycle1} is a good example of an early cycle showing fine differences in development among countries that are quite similar. Their similarity is on display in the indicator cluster maps that show them often placed in the same cluster (Fig.~\ref{4dMaps} and Fig.~\ref{Kmeans4dMaps-b}).

Note that the death filtration $\epsilon_d$ of each cycle coincides with the birth of the simplex that closes the loop. For example, consider the Northern Africa interval, $[0.85, 0.97)$, in Table~\ref{cycle2}. There are two possible internal edges: Niger~$\rightarrow$~Algeria and Mali~$\rightarrow$~Algeria, that come into existence at $\epsilon=0.94$ and $\epsilon=0.97$, respectively. All five countries make up the cycle over $[0.85,0.94)$, but it shrinks at $\epsilon=0.94$ when N~$\rightarrow$~A forms. At this point, Libya is cut off from the cycle, which persists with the four other countries until the M~$\rightarrow$~A simplex closes the loop at $\epsilon_d=0.97$.

\begin{table}[h!]
 \caption{Countries composing generating cycles and the corresponding birth and death values representing the dimension-1 homology classes of the VR complex stream built over country border graph with weights $d_I$ where $I$=\{GDP/capita (GDP), life expectancy (LE), infant mortality (IM), GNI/capita (GNI)\}. Cycles are listed in ascending order of interval birth.}
 \begin{tabular}{llp{9cm}}
 \hline\noalign{\smallskip}
 Birth & Death & Generating Countries\\
 \noalign{\smallskip}\hline\noalign{\smallskip}
0.31 & 0.52 & Hungary, Romania, Croatia, Montenegro, Serbia \\
0.46 & 0.94 & Chile, Peru, Brazil, Argentina \\
0.53 & 0.96 & Romania, Ukraine, Belarus, Poland, Hungary, Slovak Republic \\
0.54 & 0.94 & Austria, Italy, Switzerland, Germany, France \\
0.56 & 0.75 & Mali, Mauritania, Senegal, Guinea \\
0.71 & 0.85 & Congo Dem. Rep., Zambia, Tanzania, Burundi \\
0.71 & 0.81 & Kazakhstan, Turkmenistan, China, Kyrgyzstan, Uzbekistan \\
0.75 & 0.85 & China, Nepal, Bhutan, India \\
0.78 & 0.85 & Congo Dem. Rep., Uganda, Burundi, Tanzania \\
0.84 & 1.18 & Czech Rep., Germany, Austria, Slovenia, Hungary, Slovak Republic \\
0.90 & 1.38 & Congo Dem. Rep., Congo Rep., Central African Rep., Cameroon \\
0.91 & 0.96 & Syria, Turkey, Iraq, Iran \\
1.06 & 1.95 & Algeria, Mauritania, Sudan, Chad, Egypt, Niger, Mali, Libya \\
1.18 & 1.52 & Israel, Jordan, Lebanon, Syria \\
1.22 & 1.85 & Afghanistan, Turkmenistan, China, India, Tajikistan, Pakistan, Uzbekistan \\
1.24 & 1.51 & Algeria, Niger, Mauritania, Mali \\
1.26 & 1.28 & Afghanistan, Tajikistan, Turkmenistan, Uzbekistan \\
1.30 & 1.77 & Iran, Pakistan, Afghanistan, Turkmenistan \\
1.34 & 1.49 & Egypt, Israel, Jordan, Palestine \\
 \noalign{\smallskip}\hline\noalign{\smallskip}
 \end{tabular}
 \label{4d-cycles}
\end{table}

The birth of the closing simplex, i.e. the death of the cycle as a whole, indicates the overall development disparity between countries in the cycle. Compare the cycle from Chile to Bolivia with the cycle from Israel to Syria; the former interval has closing simplex at distance $d_I(\textnormal{Bolivia},\textnormal{Brazil})=0.63$, while the latter has it's close at distance $d_I(\textnormal{Israel},\textnormal{Syria})=1.16$. The greater distance between countries in the Israel cycle translates into greater developmental disparity in that region than in South America.

Even relatively short intervals can identify these local development features. The cycle of Afghanistan to Iran  (see Table~\ref{4d-cycles}) is relatively short with length $0.21$ but has one of the largest closing distances, $d_I=1.30$. This gap in development may also be visualized in Fig.~\ref{4dMaps-b} and Fig.~\ref{Kmeans4dMaps-b} in that the two countries occupy different development clusters. Hence, cycles may identify the boundaries of clusters found by persistent homology and other methods such as $K$-means.

\section{Conclusions and Further Work}
\label{sec:conclusion}
Our results show that simplistic notions of country development such as the paradigm of classifying countries as ``first'' and ``third'' world masks differences in development among countries. We find that persistent homology as a clustering algorithm does not identify two distinct clusters. PH, on the other hand, discovers fine-grained differences between countries in these categories that are hidden by Gapminder's visualization and $K$-means clustering. When comparing groups of countries we find that wealth data varies widely between groups with similar values in the health data. While countries may be below average in wealth indicators, like GDP and GNI per capita, they may have quite favorable health indicators, especially life expectancy. Bimodal paradigms conceal this fact. 

We also find geographically localized patterns that are invisible in the data consisting of just indicators in $\mathbb{R}^2$ and $\mathbb{R}^4$ by adding country border information into our analyses. First-order PH identifies cycles of development statistics among neighboring countries.  In particular, these cycles identify regions of developmental disparity, be it a subtle difference between countries as in the Brazil-Bolivia cycle or a larger gap as in the Israel-Syria cycle. Gapminder's pre-determined regions obfuscate these features as a country's membership to a region implies similarity with the other members. These cycles tell a story about development in a region that would otherwise be masked by other methods.

There are many avenues for further work with our methods. We only consider four indicators, but one may replace, add, or remove indicators as they wish to conduct studies on development or other topics. Our method allows any number of variables to be encoded either as points in a higher dimensional space or as weights in the country border graph. Gapminder hosts a bounty of indicators that may be compared in myriad combinations. 

Additionally, one may conduct a longitudinal study of persistent homology using development statistics. We use only the most recent data in our study, but there are decades worth of statistics available. Such a study would need to solve the problem of missing data. Incorporating longitudinal data would make the Betti numbers a function of both time and the filtration scale. One approach to visualize such data is the CROCKER plot discussed in \cite{Swarms}, but more techniques may become available as the theory of multi-parameter persistence is an ongoing, active area of research.

\newpage
\section*{Appendix A} \label{app:years}
\nopagebreak
\addcontentsline{toc}{section}{Appendix A}
\begin{table}[h!]
\tiny
\centering
\caption{Country and the corresponding year of the most recently-available data for each indicator, GDP per capita (GDP), Life Expectancy (LE), Infant Mortality (IM), GNI per capita (GNI).}
\begin{tabular}{ccc}
\begin{tabular}[t]{>{\raggedright\arraybackslash}p{1cm}rrrr}
 \hline\noalign{\smallskip}
Country & GDP & LE & IM & GNI \\
 \noalign{\smallskip}\hline\noalign{\smallskip}
Afghanistan & 2015 & 2016 & 2015 & 2010\\
Albania & 2015 & 2016 & 2015 & 2011\\
Algeria & 2015 & 2016 & 2015 & 2011\\
Angola & 2015 & 2016 & 2015 & 2011\\
Antigua and Barbuda & 2015 & 2016 & 2015 & 2011\\
Argentina & 2015 & 2016 & 2015 & 2011\\
Armenia & 2015 & 2016 & 2015 & 2011\\
Australia & 2015 & 2016 & 2015 & 2010\\
Austria & 2015 & 2016 & 2015 & 2011\\
Azerbaijan & 2015 & 2016 & 2015 & 2011\\
Bahamas & 2015 & 2016 & 2015 & 2010\\
Bahrain & 2015 & 2016 & 2015 & 2010\\
Bangladesh & 2015 & 2016 & 2015 & 2011\\
Barbados & 2015 & 2016 & 2015 & 2009\\
Belarus & 2015 & 2016 & 2015 & 2011\\
Belgium & 2015 & 2016 & 2015 & 2011\\
Belize & 2015 & 2016 & 2015 & 2011\\
Benin & 2015 & 2016 & 2015 & 2011\\
Bhutan & 2015 & 2016 & 2015 & 2011\\
Bolivia & 2015 & 2016 & 2015 & 2011\\
Bosnia and Herzegovina & 2015 & 2016 & 2015 & 2011\\
Botswana & 2015 & 2016 & 2015 & 2011\\
Brazil & 2015 & 2016 & 2015 & 2011\\
Brunei & 2015 & 2016 & 2015 & 2009\\
Bulgaria & 2015 & 2016 & 2015 & 2011\\
Burkina Faso & 2015 & 2016 & 2015 & 2011\\
Burundi & 2015 & 2016 & 2015 & 2011\\
Cambodia & 2015 & 2016 & 2015 & 2011\\
Cameroon & 2015 & 2016 & 2015 & 2011\\
Canada & 2015 & 2016 & 2015 & 2011\\
Cape Verde & 2015 & 2016 & 2015 & 2011\\
Central African Rep. & 2015 & 2016 & 2015 & 2011\\
Chad & 2015 & 2016 & 2015 & 2011\\
Chile & 2015 & 2016 & 2015 & 2011\\
China & 2015 & 2016 & 2015 & 2011\\
Colombia & 2015 & 2016 & 2015 & 2011\\
Comoros & 2015 & 2016 & 2015 & 2011\\
Congo Dem. Rep. & 2015 & 2016 & 2015 & 2011\\
Congo Rep. & 2015 & 2016 & 2015 & 2011\\
Costa Rica & 2015 & 2016 & 2015 & 2011\\
Cote d'Ivoire & 2015 & 2016 & 2015 & 2011\\
Croatia & 2015 & 2016 & 2015 & 2011\\
Cyprus & 2015 & 2016 & 2015 & 2010\\
Czech Rep. & 2015 & 2016 & 2015 & 2011\\
Denmark & 2015 & 2016 & 2015 & 2011\\
Djibouti & 2015 & 2016 & 2015 & 2009\\
Dominica & 2015 & 2016 & 2015 & 2011\\
Dominican Rep. & 2015 & 2016 & 2015 & 2011\\
Ecuador & 2015 & 2016 & 2015 & 2011\\
Egypt & 2015 & 2016 & 2015 & 2011\\
El Salvador & 2015 & 2016 & 2015 & 2011\\
Equatorial Guinea & 2015 & 2016 & 2015 & 2011\\
Eritrea & 2015 & 2016 & 2015 & 2011\\
Estonia & 2015 & 2016 & 2015 & 2011\\
Ethiopia & 2015 & 2016 & 2015 & 2011\\
Fiji & 2015 & 2016 & 2015 & 2011\\
Finland & 2015 & 2016 & 2015 & 2011\\
France & 2015 & 2016 & 2015 & 2011\\
Gabon & 2015 & 2016 & 2015 & 2011\\
Gambia & 2015 & 2016 & 2015 & 2011\\
\noalign{\smallskip}\hline\noalign{\smallskip}
\end{tabular}
&
\begin{tabular}[t]{>{\raggedright\arraybackslash}p{1cm}rrrr}
\hline\noalign{\smallskip}
Country & GDP & LE & IM & GNI \\
\noalign{\smallskip}\hline\noalign{\smallskip}
Georgia & 2015 & 2016 & 2015 & 2011 \\
Germany & 2015 & 2016 & 2015 & 2011 \\
Ghana & 2015 & 2016 & 2015 & 2011 \\
Greece & 2015 & 2016 & 2015 & 2011 \\
Grenada & 2015 & 2016 & 2015 & 2011 \\
Guatemala & 2015 & 2016 & 2015 & 2011 \\
Guinea & 2015 & 2016 & 2015 & 2011 \\
Guinea-Bissau & 2015 & 2016 & 2015 & 2011 \\
Guyana & 2015 & 2016 & 2015 & 2010 \\
Haiti & 2015 & 2016 & 2015 & 2011 \\
Honduras & 2015 & 2016 & 2015 & 2011 \\
Hungary & 2015 & 2016 & 2015 & 2011 \\
Iceland & 2015 & 2016 & 2015 & 2011 \\
India & 2015 & 2016 & 2015 & 2011 \\
Indonesia & 2015 & 2016 & 2015 & 2011 \\
Iran & 2015 & 2016 & 2015 & 2009 \\
Iraq & 2015 & 2016 & 2015 & 2011 \\
Ireland & 2015 & 2016 & 2015 & 2011 \\
Israel & 2015 & 2016 & 2015 & 2011 \\
Italy & 2015 & 2016 & 2015 & 2011 \\
Jamaica & 2015 & 2016 & 2015 & 2011 \\
Japan & 2015 & 2016 & 2015 & 2011 \\
Jordan & 2015 & 2016 & 2015 & 2011 \\
Kazakhstan & 2015 & 2016 & 2015 & 2011 \\
Kenya & 2015 & 2016 & 2015 & 2011 \\
Kiribati & 2015 & 2016 & 2015 & 2011 \\
Korea Rep. & 2015 & 2016 & 2015 & 2011 \\
Kuwait & 2015 & 2016 & 2015 & 2010 \\
Kyrgyzstan & 2015 & 2016 & 2015 & 2011 \\
Laos & 2015 & 2016 & 2015 & 2011 \\
Latvia & 2015 & 2016 & 2015 & 2011 \\
Lebanon & 2015 & 2016 & 2015 & 2011 \\
Lesotho & 2015 & 2016 & 2015 & 2011 \\
Liberia & 2015 & 2016 & 2015 & 2011 \\
Libya & 2015 & 2016 & 2015 & 2009 \\
Lithuania & 2015 & 2016 & 2015 & 2011 \\
Luxembourg & 2015 & 2016 & 2015 & 2011 \\
Macedonia FYR & 2015 & 2016 & 2015 & 2011 \\
Madagascar & 2015 & 2016 & 2015 & 2011 \\
Malawi & 2015 & 2016 & 2015 & 2011 \\
Malaysia & 2015 & 2016 & 2015 & 2011 \\
Maldives & 2015 & 2016 & 2015 & 2011 \\
Mali & 2015 & 2016 & 2015 & 2011 \\
Malta & 2015 & 2016 & 2015 & 2010 \\
Mauritania & 2015 & 2016 & 2015 & 2011 \\
Mauritius & 2015 & 2016 & 2015 & 2011 \\
Mexico & 2015 & 2016 & 2015 & 2011 \\
Micronesia Fed. Sts. & 2015 & 2016 & 2015 & 2011 \\
Moldova & 2015 & 2016 & 2015 & 2011 \\
Mongolia & 2015 & 2016 & 2015 & 2011 \\
Montenegro & 2015 & 2016 & 2015 & 2011 \\
Morocco & 2015 & 2016 & 2015 & 2011 \\
Mozambique & 2015 & 2016 & 2015 & 2011 \\
Namibia & 2015 & 2016 & 2015 & 2011 \\
Nepal & 2015 & 2016 & 2015 & 2011 \\
Netherlands & 2015 & 2016 & 2015 & 2011 \\
New Zealand & 2015 & 2016 & 2015 & 2010 \\
Nicaragua & 2015 & 2016 & 2015 & 2011 \\
Niger & 2015 & 2016 & 2015 & 2011 \\
Nigeria & 2015 & 2016 & 2015 & 2011 \\
Norway & 2015 & 2016 & 2015 & 2011 \\
Oman & 2015 & 2016 & 2015 & 2010 \\
Pakistan & 2015 & 2016 & 2015 & 2011 \\
\noalign{\smallskip}\hline\noalign{\smallskip}
\end{tabular}
&
\begin{tabular}[t]{>{\raggedright\arraybackslash}p{1cm}rrrr}
\hline\noalign{\smallskip}
Country & GDP & LE & IM & GNI \\
\noalign{\smallskip}\hline\noalign{\smallskip}
Palestine & 2015 & 2016 & 2015 & 2005 \\
Panama & 2015 & 2016 & 2015 & 2011 \\
Papua New Guinea & 2015 & 2016 & 2015 & 2011 \\
Paraguay & 2015 & 2016 & 2015 & 2011 \\
Peru & 2015 & 2016 & 2015 & 2011 \\
Philippines & 2015 & 2016 & 2015 & 2011 \\
Poland & 2015 & 2016 & 2015 & 2011 \\
Portugal & 2015 & 2016 & 2015 & 2011 \\
Qatar & 2015 & 2016 & 2015 & 2011 \\
Romania & 2015 & 2016 & 2015 & 2011 \\
Russia & 2015 & 2016 & 2015 & 2011 \\
Rwanda & 2015 & 2016 & 2015 & 2011 \\
Saint Lucia & 2015 & 2016 & 2015 & 2011 \\
Saint Vincent and the Grenadines & 2015 & 2016 & 2015 & 2011 \\
Samoa & 2015 & 2016 & 2015 & 2011 \\
Sao Tome and Principe & 2015 & 2016 & 2015 & 2011 \\
Saudi Arabia & 2015 & 2016 & 2015 & 2011 \\
Senegal & 2015 & 2016 & 2015 & 2011 \\
Serbia & 2015 & 2016 & 2015 & 2011 \\
Seychelles & 2015 & 2016 & 2015 & 2011 \\
Sierra Leone & 2015 & 2016 & 2015 & 2011 \\
Singapore & 2015 & 2016 & 2015 & 2011 \\
Slovak Republic & 2015 & 2016 & 2015 & 2011 \\
Slovenia & 2015 & 2016 & 2015 & 2011 \\
Solomon Islands & 2015 & 2016 & 2015 & 2011 \\
South Africa & 2015 & 2016 & 2015 & 2011 \\
Spain & 2015 & 2016 & 2015 & 2011 \\
Sri Lanka & 2015 & 2016 & 2015 & 2011 \\
Sudan & 2015 & 2016 & 2015 & 2010 \\
Suriname & 2015 & 2016 & 2015 & 2010 \\
Swaziland & 2015 & 2016 & 2015 & 2011 \\
Sweden & 2015 & 2016 & 2015 & 2011 \\
Switzerland & 2015 & 2016 & 2015 & 2011 \\
Syria & 2015 & 2016 & 2015 & 2010 \\
Tajikistan & 2015 & 2016 & 2015 & 2011 \\
Tanzania & 2015 & 2016 & 2015 & 2011 \\
Thailand & 2015 & 2016 & 2015 & 2011 \\
Timor-Leste & 2015 & 2016 & 2015 & 2010 \\
Togo & 2015 & 2016 & 2015 & 2011 \\
Tonga & 2015 & 2016 & 2015 & 2011 \\
Trinidad and Tobago & 2015 & 2016 & 2015 & 2011 \\
Tunisia & 2015 & 2016 & 2015 & 2011 \\
Turkey & 2015 & 2016 & 2015 & 2011 \\
Turkmenistan & 2015 & 2016 & 2015 & 2011 \\
Uganda & 2015 & 2016 & 2015 & 2011 \\
Ukraine & 2015 & 2016 & 2015 & 2011 \\
United Arab Emirates & 2015 & 2016 & 2015 & 2011 \\
United Kingdom & 2015 & 2016 & 2015 & 2011 \\
United States & 2015 & 2016 & 2015 & 2011 \\
Uruguay & 2015 & 2016 & 2015 & 2011 \\
Uzbekistan & 2015 & 2016 & 2015 & 2011 \\
Vanuatu & 2015 & 2016 & 2015 & 2011 \\
Venezuela & 2015 & 2016 & 2015 & 2011 \\
Vietnam & 2015 & 2016 & 2015 & 2011 \\
Yemen Rep. & 2015 & 2016 & 2015 & 2011 \\
Zambia & 2015 & 2016 & 2015 & 2011 \\
\\
 \noalign{\smallskip}\hline\noalign{\smallskip}
 \end{tabular}
 \end{tabular}
 \label{data_years}
\end{table}

\section*{Appendix B} \label{app:clusters}
\addcontentsline{toc}{section}{Appendix B}

The $K$-means clustering algorithm is an ubiquitous vector quantization method in data mining \cite{macqueen1967some}. The algorithm produces a partitioning of a point cloud of data into $K$ clusters by identifying a set of centers (or prototypes) for each of the clusters, assigning each data point to the cluster with the closest center, calculating the mean of all points in each cluster, and then updating the center of each cluster to be equal to the mean. The process is iterated until all data points are quantized to an appropriate degree of accuracy. A typical implementation of the algorithm is to randomly initialize starting centers and then to iterate for a large number of trials with the goal of minimizing the within-cluster sum of squares error.

We implement the $K$-means clustering algorithm on our datasets. We mention a couple drawbacks of this method as opposed to persistent homology: (1) The $K$-means algorithm requires a fixed number of clusters while persistent homology allows for clustering at multiple scales. Therefore, when an appropriate number of clusters is not known \emph{a priori}, specifying a set number of clusters may introduce bias. (2) Random initialization of starting centers means that a global optimum may not be achieved in the clustering, and the resulting clustering depends on this initialization, varying with different starting centers. One possible advantage of the $K$-means algorithm, however, is that elements within a cluster typically remain more similar to one another using $K$-means than the ``long'' clusters of zero-order PH.

We perform the $K$-means algorithm using random initialization of cluster centers with two numbers of clusters $K=2$ and $6$ on the $\mathbb{R}^2$ and $\mathbb{R}^4$ data sets without geographic border information. The choice of $K=2$ clusters is to consider whether the $K$-means algorithm separates countries into a first versus third world paradigm while the choice of $K=6$ coincides with the six clusters used in the Gapminder representation (see Fig. \ref{gapminder-map}) and in Section \ref{sec:H0clustering}.

\begin{figure}
\centering
\subfigure[][]{
\label{Kmeans2dMaps-a}
\includegraphics[width=0.475\linewidth]{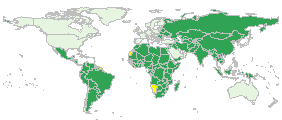}}
\hspace{8pt}
\subfigure[][]{
\label{Kmeans2dMaps-b}
\includegraphics[width=0.475\linewidth]{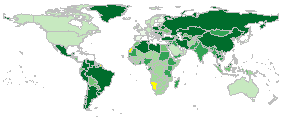}} \\
\caption{World map depicting clusters found using $K$-means of $\mathbb{R}^2$ data:
\subref{Kmeans2dMaps-a} $K=2$ with cluster sizes 59 and 135;
\subref{Kmeans2dMaps-b} $K=6$ with cluster sizes 21, 25, 31, 32, 36, and 49.
Shade corresponds to cluster size, where darker is larger. Yellow denotes countries missing from the data set because not all indicators are available.}
\label{Kmeans2dMaps}
\end{figure}

\begin{figure}
\centering
\subfigure[][]{
\label{Kmeans4dMaps-a}
\includegraphics[width=0.475\linewidth]{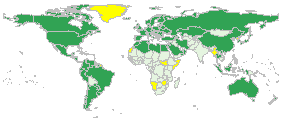}}
\hspace{8pt}
\subfigure[][]{
\label{Kmeans4dMaps-b}
\includegraphics[width=0.475\linewidth]{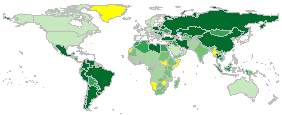}} \\
\caption{World map depicting clusters found using $K$-means of $\mathbb{R}^4$ data:
\subref{Kmeans4dMaps-a} $K=2$ with cluster sizes 71 and 108;
\subref{Kmeans4dMaps-b} $K=6$ with cluster sizes 19, 20, 21, 37, 40, and 42.
Shade corresponds to cluster size, where darker is larger. Yellow denotes countries missing from the data set because not all indicators are available.}
\label{Kmeans4dMaps}
\end{figure}

In Fig.~\ref{Kmeans2dMaps-a}, we observe that the $K=2$ clustering of the $\mathbb{R}^2$ data appears to follow what some may view as a first versus third world paradigm, grouping wealthier countries together. However, once more indicators have been added, this distinction starts to break down as the additional indicators reveal a more nuanced notion of similarity, see Fig.~\ref{Kmeans4dMaps-a}. We observe that in Fig.~\ref{Kmeans2dMaps-b} and \ref{Kmeans4dMaps-b}, when a more fine-grained approach is used to split the countries into more clusters, the clusters do not split along this paradigm or traditional continental divisions.

\bibliographystyle{plain}
\bibliography{bibliography}

\end{document}